\documentclass{article}
\usepackage{amsmath}
\usepackage{amssymb}
  \usepackage{paralist}
  \usepackage{graphics} 
  \usepackage{epsfig} 
\usepackage{graphicx}  \usepackage{epstopdf}
 \usepackage[colorlinks=true]{hyperref}
\hypersetup{urlcolor=blue, citecolor=red}
\usepackage{mathrsfs}
\usepackage{booktabs}

\newtheorem{theorem}{Theorem}[section]

\newtheorem{lemma}[theorem]{Lemma}
\newtheorem{proposition}[theorem]{Proposition}

\newtheorem{definition}[theorem]{Definition}
\newtheorem{remark}{Remark}

\title{A Crank-Nicolson type minimization scheme \\ for a hyperbolic free boundary problem}

\date{}
\begin{document}
\maketitle

\centerline{\scshape Yoshiho Akagawa}
\medskip
{\footnotesize
 \centerline{National Institute of Technology, Gifu College}
   \centerline{Motosu, Gifu, 501-0495 Japan}
   \centerline{akagawa@gifu-nct.ac.jp}
} 

\medskip

\centerline{\scshape Elliott Ginder}
\medskip
{\footnotesize
 \centerline{Meiji Institute for Advanced Study of Mathematical Sciences, Meiji University}
   \centerline{Nakanoku, Tokyo, 164-8525 Japan}
      \centerline{eginder@meiji.ac.jp}
} 

\medskip

\centerline{\scshape Syota Koide}
\medskip
{\footnotesize
 \centerline{Graduate School of Natural Science and Technology, Kanazawa University}
   \centerline{Kanazawa, Ishikawa, 920-1192 Japan}
      \centerline{math.koide@gmail.com}
} 

\medskip

\centerline{\scshape Seiro Omata}
\medskip
{\footnotesize
 \centerline{Institute of Science and Engineering, Kanazawa University}
   \centerline{Kanazawa, Ishikawa, 920-1192 Japan}
      \centerline{omata@se.kanazawa-u.ac.jp}
} 

\medskip

\centerline{\scshape Karel Svadlenka{\scriptsize (corresponding author)}}
\medskip
{\footnotesize
 \centerline{Graduate School of Science, Kyoto University}
   \centerline{Sakyoku, Kyoto, 606-8502 Japan}
      \centerline{karel@math.kyoto-u.ac.jp}
}

\bigskip

\begin{abstract}
We consider a hyperbolic free boundary problem by means of minimizing 
time discretized functionals of Crank-Nicolson type. The feature of this functional is that 
it enjoys energy conservation in the absence of free boundaries, 
which is an essential property for numerical calculations.
The existence and regularity of minimizers is shown and an energy estimate is derived.
These results are then used to show the existence of a weak solution to the free boundary problem in 
the 1-dimensional setting. 
\end{abstract}

\section{Introduction}

In this paper we treat a variational problem related to 
the following hyperbolic free boundary problem:

\noindent{\bf{Problem (1).}}
Find $u : \Omega \times [0, T) \to \mathbb{R}$ such that
\begin{equation}
\begin{cases}
\chi_{\overline{\{u>0\}}}\, u_{tt} &= \Delta u \quad \hbox{in} \,\,\Omega \times (0, T), \\
u(x, 0) &= u_0(x) \quad \hbox{in} \,\, \Omega, \\
u_t(x, 0) &= v_0(x) \quad \hbox{in} \,\, \Omega,\\
\end{cases}
\label{prob}
\end{equation}
under suitable boundary conditions, 
where $\Omega \subset \mathbb{R}^{N}$ is a bounded Lipschitz domain, $T>0$ is the final time, $u_0$  denotes the initial condition, $v_0$ is the initial velocity, and $\{u >0\}$ is the set $\{(x, t) \in \Omega \times 
(0, T): u(x, t)>0\}$.

We observe that $u$ satisfies the wave equation where $u>0$, and that $u$ is harmonic when $u<0.$ 
It can be formally shown that, in the energy-preserving regime, solutions fulfill
the following free boundary condition: 
\begin{equation}
|\nabla u|^2 -u^2_t = 0 \quad \hbox{on}\,\,\Omega \times (0, T) \cap \partial \{u>0\}. 
\label{eq:fbc0}
\end{equation}

This kind of problem is a natural prototype for explaining phenomena  involving oscillations in the presence of an obstacle, e.g., an elastic string hitting a desk or soap bubbles moving atop water. 
However, due to the lack of rigorous mathematical results for hyperbolic free boundary problems, numerical studies of this problem are very limited. 
On the other hand, interesting results have been already obtained for first-order hyperbolic free boundary problems, see, e.g., \cite{[CCF]} for the analysis of a model of tumor growth, or \cite{[IL]} for the analysis of wave-structure interactions.

In this paper, we propose a numerical scheme with good energy conservation properties and provide its theoretical background at least in the case of space dimension 1.
We remark that similar types of problems have been treated in \cite{[KO]}, \cite{[YOSO]}, \cite{[Kikuchi]}
\cite{[G-S]},
and that the recent paper \cite{[O4]} has established a precise mathematical formulation. 
These papers revealed that the discrete Morse flow (also known as minimizing movements), a variational method based on time-discretized functionals, is an effective tool not only for problems of elliptic and parabolic type but also in the hyperbolic setting.

We now briefly review those previous results.
Kikuchi and Omata \cite{[KO]} studied the problem in the one-dimensional domain 
$\Omega = (0, \infty)$. 
They showed the existence and uniqueness of a strong solution $u \in C^2 (\Omega \times (0, 
\infty) \cap \{u>0\})$, the regularity of its free boundary $\partial \{u>0\}$ and the well-posedness of the
problem under suitable compatibility conditions. 
Yoshiuchi et al. \cite{[YOSO]} addressed a similar problem to Problem \ref{prob}
 including a damping term 
$\alpha u_t$. Using the discrete Morse flow, they derived an energy estimate for approximate solutions, and provided numerical results. 
The following problem, stated here without initial and boundary conditions, 
has been treated by Kikuchi in \cite{[Kikuchi]}: 
\begin{equation*}
\begin{cases}
& u_{tt} -u_{xx} \geq  0\quad \hbox{in} \,\, (0,1) \times (0, \infty) \cap \{u>0\}, \\
&\hbox{spt}\,(u_{tt}-u_{xx}) \subset \{u=0\}, \\
&u(x, t) \geq 0  \quad \hbox{$L^2$-a.e.}. 
\end{cases}
\end{equation*}
He constructed a weak solution to this problem using a minimizing method in the 
spirit of the discrete Morse flow. Moreover, two of the present authors, 
Ginder and Svadlenka, dealt with a hyperbolic free boundary problem under a  
volume constraint in \cite{[G-S]}. 
They constructed a weak solution in the one dimensional setting, again using the discrete 
Morse flow. 
The analysis of hyperbolic obstacle problems based on this variational framework 
of discrete Morse flow was significantly extended to nonlocal problems
(with fractional Laplacian) and to semilinear problems in a vector-valued setting by the group
of Bonafini, Novaga and Orlandi in the papers \cite{[BNO], [BLNO]}.

The technique used in this paper is similar to the discrete Morse flow used in the 
above papers, but it has some distinct differences. 
We now explain our approach, together with the organization of the paper. 
The first step is to consider the minimization, within the set
$$
\mathcal{K}:= \{u \in H^1 (\Omega) \,;\,u = u_0 \,\hbox{on}\,\partial \Omega\}, 
$$
of the following time-discretized, Crank-Nicolson type functional: 
\begin{equation}
\label{CN-functional}
J_m(u) := \int_{\Omega \cap \mathcal{S}_m(u)}
\frac{|u -2u_{m-1} +u_{m-2}|^2}{2 h^2} \,dx + \frac{1}{4} \int_{\Omega} 
|\nabla u + \nabla u_{m-2}|^2 \,dx 
\end{equation}
where $\mathcal{S}_m(u)$ denotes the set $\{u>0\} \cup \{u_{m-1}>0\} \cup \{u_{m-2} >0\}$. 
In the above, $u_m$ 
represents an approximation of the solution to the original problem at a fixed time $t_m = mh$, 
where $m\in{\mathbb N}$, and $h>0$ denotes the time step, obtained by dividing the time 
interval $(0,T)$ into $M$ equal parts.
The difference with the previous approaches (e.g., \cite{[YOSO]}, \cite{[G-S]}) appears in the 
gradient term. As will be shown, our gradient term yields the energy conservation property 
even in the time 
discretized setting when there is no free boundary (i.e., when the restriction to the set 
$\mathcal{S}_m(u)$ is omitted). 
This is a significant improvement to the previous approaches, where energy conservation
does not hold and leads to inaccurate numerical solutions.  
We believe that combination of the new energy-preserving scheme with the more advanced 
analytical techniques developed in \cite{[BNO], [BLNO]} will lead to a powerful framework.
The details of the new variational scheme are presented in Section \ref{sec:EC}. 

The sequence $\{u_m\}$ is constructed using the initial conditions 
$u_0 \in \mathcal{K}$ and $v_0 \in H^1_0(\Omega)$. 
In particular, a forward differencing is used to define $u_1:= u_0 + h v_0 \in {\mathcal K}$. 
Then, for any integer $m \geq 2$, a minimizer $\widetilde{u}_m$ of $J_m$ exists and has
the subsolution property (Theorem \ref{exi_mini}, Proposition \ref{subsolution} in 
Section \ref{sec:minimizing method}), 
which can be shown by a standard argument, such as in \cite{[YOSO]}. 
We then define the cut-off minimizer $u_m := \max \{\widetilde{u}_m, 0 \}$ and find that the 
minimizers $u_m$ satisfy an energy estimate (Theorem \ref{enerest}).
In particular, for any integer $k \geq 1$, we have 
$$
\Bigl \| \frac{u_k - u_{k-1}}{h}\Bigr\|_{L^2(\Omega)}^2 + \frac{1}{2} \|\nabla u_k\|
_{L^2(\Omega)}^2 \leq C(u_0,v_0, \Omega),
$$
where $C$ is independent of $k$ and $h$.
The regularity of $u_m$ can then be obtained and this leads to the first variation 
formula (Proposition \ref{first_vari_form}).  
The next step is to define approximate weak solutions of the free boundary problem, as well as 
a weak solution of Problem \ref{prob}. This is done in Section \ref{sec:weak solution}. 
In the case of a 1-dimensional domain, the energy estimate and the embedding 
theorem allow us to pass the time step size $h$ to zero in the equation of approximate 
weak solutions to obtain a weak solution of Problem \ref{prob}. 

In Section \ref{sec:weak solution}, we also formally discuss the role of free boundary condition 
on a more general setting for this problem, namely 
a hyperbolic free boundary problem with an adhesion term. 
When energy is not conserved, we cannot calculate the first variation in the usual way. 
In this case, we adopt a formal approach and 
derive the equation from a measure theoretic point of view. 
We conclude by presenting numerical results and their analysis in Section \ref{sec:numerics}, 
emphasizing the wide applicability of the proposed numerical method.

\section{Energy conservation}
\label{sec:EC}
In this section, we derive the energy preserving property of the 
functional $J_m$ in (\ref{CN-functional}), which has not been achieved in previous research.
To this end, let us consider the following modified functional:
\begin{equation}
\label{eq:ECfunctional}
I_m(u) := \int_{\Omega}
\frac{|u -2u_{m-1} +u_{m-2}|^2}{2 h^2} \,dx + \frac{1}{4} \int_{\Omega} 
|\nabla u + \nabla u_{m-2}|^2 \,dx, 
\end{equation}
on the set $\mathcal{K}:= \{u \in H^1 (\Omega) \,;\,u = u_0 \,\hbox{on}\,\partial \Omega\}$. 
This functional can be regarded as the no-free boundary version of $J_m$. We
note that a unique minimizer exists for each $I_m$ whenever $I_m(u_0)<\infty$ 
since the functional is convex and lower semicontinuous. 

\begin{theorem}[Energy conservation]
Minimizers $u_k$ of $I_k$ conserve the energy
\begin{equation}
E_k := \Bigl\|\frac{u_k -u_{k-1}}{h} \Bigr\|^2_{L^2(\Omega)} + \frac{1}{2}
 \Bigl( \|\nabla u_k\|^2_{L^2(\Omega)} +\|\nabla u_{k-1}\|^2_{L^2(\Omega)} \Bigr), 
\label{eq:ek}
\end{equation}
in the sense that $E_k$ is independent of $k \geq 1$. 
\end{theorem}

{\emph{Proof. }}
For $m\geq 2$, the function $(1-\lambda) u_m + \lambda u_{m-2} = u_m + \lambda(u_{m-2}-u_{m})$
is admissible for every $\lambda \in [0, 1]$, which justifies 
$$
\left. \frac{d}{d \lambda} I_m(u_m + \lambda(u_{m-2}-u_{m})) \right|_{\lambda =0}=0.
$$
Computing this derivative, we have 
\begin{align*}
0 &= \frac{d}{d\lambda} \int_{\Omega} 
\Bigl[\frac{|u_{m}+\lambda(u_{m-2}-u_{m})-2u_{m-1}+u_{m-2}|^2}{2h^2} \\
&\left. \qquad\qquad+ \frac{1}{4}|\nabla (u_m + \lambda(u_{m-2}-u_m))+\nabla u_{m-2}|^2 \Bigr] \right|
_{\lambda = 0}\,dx \\
&=\int_{\Omega}\Bigl[\frac{(u_{m-2}-u_{m}) (u_{m}-2u_{m-1}+u_{m-2})}{h^2} \\
& \qquad \qquad +\frac{1}{2}\nabla(u_{m-2}-u_m) \cdot \nabla(u_m + u_{m-2}) \Bigr]\, dx \\
& =\int_{\Omega}\Bigl[\frac{(u_{m-1}-u_{m-2})^2 - (u_{m}-u_{m-1})^2}{h^2}
+\frac{1}{2}|\nabla u_{m-2}|^2-  \frac{1}{2} |\nabla u_m|^2 \Bigr]\, dx.
\end{align*}
\noindent
Summing over $m=2, ..., k$, we arrive at 
\begin{multline*}
\int_{\Omega}
\Bigl[ \frac{1}{h^2}(u_{1}-u_{0})^2 -  \frac{1}{h^2}(u_{k}-u_{k-1})^2  \\
+ \frac{1}{2}  | \nabla u_{0}|^2+\frac{1}{2} | \nabla u_{1}|^2
 -\frac{1}{2} | \nabla u_{k-1}|^2 - \frac{1}{2} | \nabla  u_{k} |^2 
\Bigr] \,dx = 0,
\end{multline*}
which means 
\begin{eqnarray*}
E_1 &=& \int_{\Omega}
\Bigl[ \frac{1}{h^2}(u_{1}-u_{0})^2 +
\frac{1}{2}  | \nabla u_{0}|^2+ \frac{1}{2}  | \nabla u_{1}|^2\Bigr] \,dx \\
&=& \int_{\Omega} \Bigl[ \frac{1}{h^2}(u_{k}-u_{k-1})^2 
+ \frac{1}{2}  | \nabla u_{k-1}|^2+\frac{1}{2} | \nabla  u_{k} |^2 \Bigr] \,dx \\
&=& E_k, 
\end{eqnarray*}
and the proof is complete.

\section{The minimizing method}
\label{sec:minimizing method}
For any integer $m \geq 2$, we introduce the following functional:
\begin{equation}
J_m(u)
= \int_{\Omega \cap \mathcal{S}_m(u)}
\frac{|u -2u_{m-1} +u_{m-2}|^2}{2 h^2} \,dx 
+ \frac{1}{4} \int_{\Omega} 
|\nabla u + \nabla u_{m-2}|^2 \,dx, 
\label{eq:2.1}
\end{equation}
where $\mathcal{S}_m(u) := \{u>0\} \cup \{u_{m-1}>0\} \cup \{u_{m-2} >0\}$.

We determine a sequence of functions $\{ u_{m} \}$ iteratively
by taking $u_0 \in {\mathcal K}$ and $u_1 = u_0 + h v_0 \in {\mathcal K}$, 
defining $\widetilde{u}_m$ as a minimizer of $J_m$
in $\mathcal K$, and setting $u_{m} := \max\{\widetilde{u}_m, 0\}$.

We now study the existence and regularity of minimizers
which guarantees the possibility of applying the first variation formula
to $J_m$.

\begin{theorem}[Existence]
\label{exi_mini}
If $J_m(u_0)<\infty$, then there exists a minimizer $\widetilde{u}_{m} \in {\mathcal K}$ of
the functional $J_m$.
\end{theorem}

{\emph{Proof. }}
Given $u_{m-1}, u_{m-2}$, we show the existence of $\widetilde{u}_m$. 
Since the infimum of $J_m$ is non-negative, we have only to show
the lower semi-continuity of $J_m$. Take any minimizing sequence $\{u^j\} \subset \mathcal{K}$ 
such that $J_m(u^j) \to \inf_{u \in \mathcal{K}} J_m(u) $ as $j \to \infty$. 
Since the sequence $\{u^j -u_0\} \subset H^1_0(\Omega)$ 
is bounded in $H^1(\Omega)$,
there exist $\widetilde{u} \in H^1_0 (\Omega)$ and $\gamma 
\in L^p(\Omega)$, $p \in [1, \infty)$ such that, up to extracting a subsequence, 
\begin{align}
\label{mini_proof1}
u^j -u_0 &\to \,\widetilde{u} \quad \hbox{strongly in} \,\,L^2(\Omega), \nonumber \\
\nabla (u^j -u_0) &\rightharpoonup \, \nabla \widetilde{u} \quad  \hbox{weakly in} \,\,L^2(\Omega), \\
\chi_{\mathcal{S}_m(u^j)} &\rightharpoonup \,\gamma \quad  \hbox{weakly  in} 
\,\,L^p(\Omega) \nonumber. 
\end{align}
Since $0 \leq \gamma \leq 1$ a.e. on $\Omega$, and $\gamma=1$ a.e. on $\{u>0\}$, where
$u := \widetilde{u} + u_0 \in \mathcal{K}$, 

\begin{align*}
J_m(u) &= \int_{\Omega}
\frac{|u -2u_{m-1} +u_{m-2}|^2}{2 h^2} \chi_{\mathcal{S}_m(u)} \,dx 
+ \frac{1}{4} \int_{\Omega} 
|\nabla u + \nabla u_{m-2}|^2 \,dx \\
&\leq \int_{\Omega}
\frac{|u -2u_{m-1} +u_{m-2}|^2}{2 h^2} \gamma \,dx 
+ \frac{1}{4} \int_{\Omega} 
|\nabla u + \nabla u_{m-2}|^2 \,dx \\
& \leq \liminf_{j \to \infty} J_m(u^j), 
\end{align*}
where the second inequality follows from (\ref{mini_proof1}).

The minimizers of $J_m$ have the following subsolution property. 

\begin{proposition}[Subsolution]
\label{subsolution}
Any minimizer $u$ of $J_m$
satisfies the following inequality
for arbitrary nonnegative $\zeta \in H_0^1(\Omega)$:
\begin{equation}
\label{subsol1}
\int_{\Omega \cap \mathcal{S}_m(u)} 
\frac{u-2u_{m-1}+u_{m-2}}{h^2} \zeta \, dx 
+ \int_{\Omega} \nabla \frac{u+u_{m-2}}{2} \cdot \nabla \zeta \, dx \leq 0.
\end{equation}
\end{proposition}
{\emph{Proof. }}
Fixing $\zeta \in C_0^{\infty}(\Omega)$ with $\zeta \geq 0$, and $\varepsilon>0$, 
we have 
\begin{align}
\label{subpro_proof1}
0 &\leq J_m(u -\varepsilon \zeta) - J_m(u) \quad (\hbox{by the minimality of $u$}) \nonumber \\
&= \int_{\Omega}
\frac{|(u -\varepsilon \zeta) -2u_{m-1} +u_{m-2}|^2}{2 h^2} \chi_{\mathcal{S}_m(u -\varepsilon \zeta)} 
\,dx + \frac{1}{4} \int_{\Omega} 
|\nabla (u-\varepsilon \zeta) + \nabla u_{m-2}|^2 \,dx \nonumber  \\
&\qquad  - \Bigl( \int_{\Omega}
\frac{|u -2u_{m-1} +u_{m-2}|^2}{2 h^2} \chi_{\mathcal{S}_m(u)} 
\,dx + \frac{1}{4} \int_{\Omega} 
|\nabla u + \nabla u_{m-2}|^2 \,dx  \Bigr). 
\end{align}
Noting that,
\begin{align*}
&\chi_{\mathcal{S}_m(u -\varepsilon \zeta)} - \chi_{\mathcal{S}_m(u)} \leq 0, \\
&|(u-\varepsilon \zeta) -2 u_{m-1} +u_{m-2}|^2
-|u -2 u_{m-1} +u_{m-2}|^2 \\
&\qquad \qquad \qquad  \qquad \qquad   \qquad \qquad   \qquad \qquad = -2 \varepsilon \zeta (u-2 u_{m-1} +u_{m-2}) + \varepsilon^2 \zeta^2 , \\
 &|\nabla (u - \varepsilon \zeta) + \nabla u_{m-2}|^2 - |\nabla u + \nabla u_{m-2}|^2 = 
-2 \varepsilon (\nabla u + \nabla u_{m-2}) \cdot \nabla \zeta + \varepsilon^2 |\nabla \zeta|^2, 
\end{align*}
we continue the estimate as
\begin{align*}
(\ref{subpro_proof1}) &\leq 
\int_{\Omega} \{ -2 \varepsilon \zeta (u-2 u_{m-1} +u_{m-2}) + \varepsilon^2 \zeta^2 \} 
\times \frac{1}{2h^2} \chi_{\mathcal{S}_m(u)} \,dx  \\
&\qquad 
+ \frac{1}{4} \int_{\Omega} \{ -2 \varepsilon (\nabla u + \nabla u_{m-2}) \cdot \nabla \zeta + \varepsilon^2 |\nabla \zeta|^2 \}\,dx. 
\end{align*}
Dividing by $\varepsilon$, letting $\varepsilon$ decrease to zero from above, and applying a density 
argument concludes the proof.

We now derive an energy estimate satisfied by the minimizers of $J_m$. 

\begin{theorem}[Energy estimate]
\label{enerest}
For any integer $k \geq 1$, we have
\begin{equation}
\label{EE0}
\Bigl \| \frac{u_k - u_{k-1}}{h}\Bigr\|_{L^2(\Omega)}^2 + \frac{1}{2} \|\nabla u_k\|
_{L^2(\Omega)}^2 \leq 
\|v_0 \|_{L^2(\Omega)}^2 + \frac{1}{2} \|\nabla u_0\|
_{L^2(\Omega)}^2 + \frac{1}{2} \|\nabla u_1\|_{L^2(\Omega)}^2.
\end{equation}
\end{theorem}
{\emph{Proof. }}
Since the function $(1-\lambda)\widetilde{u}_m + \lambda u_{m-2} = \widetilde{u}_{m} + \lambda
(u_{m-2} -\widetilde{u}_m)$ belongs to $\mathcal{K}$ for any $\lambda \in [0, 1]$, 
by the minimality property, we have
$J_m (\widetilde{u}_{m}) \le J_m(\widetilde{u}_{m} + \lambda(u_{m-2} -\widetilde{u}_m))$, and 
thus,
\begin{equation}
\label{EE1}
\lim_{\lambda \rightarrow 0+} {\frac{1}{\lambda}} \Bigl( J_m (\widetilde{u}_{m}+\lambda(u_{m-2}-\widetilde{u}_{m}))
-J_m (\widetilde{u}_{m}) \Bigr) \ge 0.
\end{equation}
Let $A_m$ denote the set 
$$
A_m := \Omega \cap (\{\widetilde{u}_{m}>0\}\cup\{u_{m-1}>0\}\cup\{u_{m-2}>0\}). 
$$
We investigate the behavior of the individual terms in $(\ref{EE1})$. For the
gradient term we get
\begin{eqnarray}
\label{grad_es}
&& \lim_{\lambda \rightarrow 0+} \frac{1}{4 \lambda} \left( | \nabla
( \widetilde{u}_{m} + \lambda (u_{m-2}- \widetilde{u}_{m})+u_{m-2})|^2 - | \nabla (\widetilde{u}_{m}+u_{m-2})
|^2 \right)  \nonumber \\
&& \qquad \qquad = \frac{1}{2} \nabla (\widetilde{u}_{m}+ u_{m-2}) \cdot \nabla
 (u_{m-2}- \widetilde{u}_{m}) \, dx \nonumber \\
&& \qquad \qquad = \frac{1}{2}  | \nabla u_{m-2}|^2  - \frac{1}{2} | \nabla \widetilde{u}_{m} |^2 \nonumber \\
&& \qquad \qquad \leq \frac{1}{2}  | \nabla u_{m-2}|^2
 - \frac{1}{2} | \nabla  u_{m} |^2.
\end{eqnarray}
For the time-discretized term, 
taking into account that the set 
$$
B_m(\lambda):= \{\widetilde  u_{m}+\lambda (u_{m-2}-\widetilde  u_{m})>0\} 
\cup\{u_{m-1}>0\}\cup\{u_{m-2}>0\}
$$
is contained in the set $A_m$, 
we find that
\begin{align*}
&\frac{1}{2h^{2}}\int_{\Omega}
\Bigl( |\widetilde{u}_{m}+\lambda(u_{m-2}-\widetilde{u}_{m})-2u_{m-1}+u_{m-2}|^{2} 
\chi_{B_m(\lambda)} \\
&\qquad - |\widetilde{u}_{m}-2u_{m-1}+u_{m-2}|^{2}
\chi_{A_m} \Bigr) \, dx \\
&\le
\frac{1}{2h^{2}}\int_{\Omega}
\Bigl( |\widetilde{u}_{m}+\lambda(u_{m-2}-\widetilde{u}_{m})-2u_{m-1}+u_{m-2}|^{2} \\
&\qquad - |\widetilde{u}_{m}-2u_{m-1}+u_{m-2}|^{2} \Bigr) \chi_{A_m} \, dx. 
\end{align*}
\noindent
Then we have
\begin{align}
\label{time_dis_estimate}
&\lim_{\lambda \to 0+} \frac{1}{2h^2 \lambda} 
\int_{\Omega}
\Bigl( |\widetilde{u}_{m}+\lambda(u_{m-2}-\widetilde{u}_{m})-2u_{m-1}+u_{m-2}|^{2}\nonumber \\
&\qquad - |\widetilde{u}_{m}-2u_{m-1}+u_{m-2}|^{2} \Bigr) \chi_{A_m}  \, dx \nonumber \\
&=\lim_{\lambda \to 0+} \frac{1}{2h^2}  \int_{A_m}
(u_{m-2}-\widetilde{u}_{m}) ( 2\widetilde{u}_{m}+\lambda(u_{m-2}-\widetilde{u}_{m})-4u_{m-1}+2u_{m-2}) \, dx  \\
& =\frac{1}{h^2}  \int_{A_m}
(u_{m-2}-\widetilde{u}_{m}) ( \widetilde{u}_{m}-2u_{m-1}+u_{m-2}) \, dx  \nonumber\\
&=\frac{1}{h^2}  \int_{A_m}
[(u_{m-1}-u_{m-2})^2 -  (u_{m-1}-\widetilde{u}_{m})^2]  \,dx. \nonumber
\end{align}
Now, $\int_{A_m} (u_{m-1}-u_{m-2})^2 \,dx \leq \int_{\Omega} (u_{m-1}-
u_{m-2})^2 \,dx$ since the integrand is non-negative. Moreover, $u_m = \max \{\widetilde{u}_m, 0
\}$ and $u_{m-1} \geq 0$ imply $(u_{m-1} -\widetilde{u}_{m})^2 \geq (u_{m-1} -u_{m})^2$, therefore
$$
-\int_{A_m} (u_{m-1} -\widetilde{u}_{m})^2\,dx \leq 
-\int_{A_m} (u_{m-1} -u_{m})^2\,dx.
$$
Noting that, outside of $A_m$, both $u_m$ and $u_{m-1}$ vanish, we get
$$
-\int_{A_m} (u_{m-1} -u_{m})^2\,dx
 =-\int_{\Omega} (u_{m-1} -u_{m})^2\,dx.
$$
Returning to (\ref{time_dis_estimate}), we get the estimate for the time discretized term:
$$
\hbox{the right hand side of}\,\,(\ref{time_dis_estimate}) \leq \frac{1}{h^2}  \int_{\Omega}
[(u_{m-1}-u_{m-2})^2 -  (u_{m-1}- u_{m})^2] dx.
$$
Combining this result and the gradient term estimate (\ref{grad_es}), we obtain
$$
 \int_{\Omega}
\Bigl[ \frac{1}{h^2}(u_{m-1}-u_{m-2})^2 -  \frac{1}{h^2}(u_{m-1}- u_{m})^2 +
\frac{1}{2}  | \nabla u_{m-2}|^2
 - \frac{1}{2} | \nabla  u_{m} |^2
\Bigr] \,dx \geq 0.
$$
Summing over $m=2, ..., k$, we arrive at 
\begin{multline*}
\int_{\Omega}
\Bigl[ \frac{1}{h^2}(u_{1}-u_{0})^2 -  \frac{1}{h^2}(u_{k}- u_{k-1})^2 \\
+ \frac{1}{2}  | \nabla u_{0}|^2+ \frac{1}{2}  | \nabla u_{1}|^2
 - \frac{1}{2}  | \nabla u_{k-1}|^2 -\frac{1}{2} | \nabla  u_{k} |^2
\Bigr] \,dx \geq 0,
\end{multline*}
which, after omitting the term $|\nabla u_{k-1}|^2 \geq 0,$ yields the desired estimate.

The following theorem is obtained by a standard argument from
elliptic regularity theory. For the sake of completeness, we shall briefly demonstrate it.  

\begin{theorem}[Regularity]
\label{regular}
Assume, in addition, that  $u_0, u_1$ belong to $L^\infty(\Omega) \cap 
C^{0, \alpha_0}_{\mathrm{loc}}(\Omega)$ 
for some $\alpha_0 \in (0, 1)$, where $u_1 := u_0 + h v_0$, and $u_0$ are non-negative. 
For every $\widetilde{\Omega} \subset \subset \Omega$, there exists
a positive constant $\alpha \in (0, 1)$ independent of $m$, 
such that the minimizers $\widetilde{u}_{m}$ belong to $C^{0, \alpha}(\tilde{\Omega})$.
\end{theorem}

To prove this, we prepare two lemmas. 
\begin{lemma}
\label{bounded}
$\widetilde{u}_m \in L^\infty(\Omega)$ 
for every $m \geq 2$. 
\end{lemma}

{\emph{Proof. }}
We use mathematical induction for $m \geq 2$. 
For $m = 2$, setting $\psi_\delta (u):= u - \delta (u + u_0 -k)^+ \in \mathcal{K}$, 
where $u := \widetilde{u}_2$, $(u+u_0 -k)^+ :=\max \{u+u_0 -k, 0\} $, 
$\delta >0$, $k \geq \max\{2 \max_{\partial \Omega} u_0, 1\}$, 
we calculate the quantity $J_2(\psi_\delta(u))-J_2(u)$, which is non-negative by the 
minimality of $u$. Noting that $S_m(\psi_{\delta}(u)) \subset S_m(u)$, we have
\begin{align*}
0 &\leq J_2(\psi_\delta(u))-J_2(u) \\
&\leq \int_{\Omega} \Bigl(\frac{|\psi_{\delta}(u) - 2 u_{1} + u_{0}|^2}{2 h^2}
-  \frac{|u - 2 u_{1} + u_{0}|^2}{2 h^2} \Bigr) \chi_{\mathcal{S}_2 (u)}\,dx \\
&\qquad + \frac{1}{4} \int_{\Omega} \Bigl( |\nabla \psi_{\delta}(u) + \nabla u_0|^2 - 
|\nabla u + \nabla u_0|^2 \Bigr) \, dx. 
\end{align*}
Dividing by $\delta$, letting $\delta \to 0+$, and setting $A_k := \{u+u_0 > k\}$, 
we get
\begin{align*}
0  &\leq - \int_{A_{k} \cap \mathcal{S}_2(u)} \frac{u - 2 u_1 + u_0}{h^2} (u+ u_0 -k) \,dx 
- \frac{1}{2} \int_{A_k} |\nabla u + \nabla u_0|^2 \,dx \\
& \leq \int_{A_{k} \cap \mathcal{S}_2(u)} \frac{2 u_1}{h^2}(u+ u_0 -k) \,dx
- \frac{1}{2} \int_{A_k} |\nabla u + \nabla u_0|^2 \,dx \\
& \leq \frac{C}{h^2} \Bigl( \frac{1}{2} \int_{A_{k}} (u+ u_0 -k)^2 \,dx + \frac{1}{2} |A_k| \Bigr) 
- \frac{1}{2} \int_{A_k} |\nabla u + \nabla u_0|^2 \,dx, 
\end{align*}
where we have used Young's inequality at the last line. Since $k \geq 1$, we get 
$$
\int_{A_k} |\nabla u + \nabla u_0|^2 \,dx \leq C 
\Bigl( \int_{A_{k}} (u+ u_0 -k)^2 \,dx + k^2 |A_k|  \Bigr). 
$$
Therefore, by \cite[Theorem 2.5.1]{[L-U]}, 
we find that $u + u_0 \in L^\infty(\Omega)$, and hence $u = \widetilde{u}_2 \in L^\infty(\Omega)$.  

Next, we assume that $\widetilde{u}_k \in L^\infty(\Omega)$ for all $k=2, ..., m-1$. 
Since $u_k = \max \{\widetilde{u}_k, 0 \}$ belongs to 
$L^\infty(\Omega)$ for all $k=2, ..., m-1$, 
by repeating the above argument with $\widetilde{u}_2, u_1, u_0$ replaced by 
$\widetilde{u}_m, u_{m-1}, u_{m-2}$, respectively, 
we get $\widetilde{u}_m + u_{m-2} \in L^\infty(\Omega)$. Therefore, $\widetilde{u}_m \in L^\infty(\Omega)$.

The previous result implies that there exists  $\mu >0$, 
which depends only on $\Omega, u_0, u_1, h$ but not on $m$, 
such that  $\sup_\Omega |\widetilde{u}_m| \leq \mu$. 

\begin{lemma}
\label{Cacciopoli_ineq}
Fix $d>0$. There exists $\gamma = \gamma(\Omega, \mu, d, h) >0$ such that for $U = \pm (\widetilde{u}_m + u_{m-2})$, 
  $$
  \int_{A_{k, r-\sigma r}} |\nabla U|^2 \,dx \leq \gamma \Bigl[\frac{1}{(\sigma r)^2}\sup_{B_r} 
  (U-k)^2 + 1 \Bigr] |A_{k, r}|
  $$
  for all $\sigma \in (0, 1)$, $B_r \subset \Omega$, and $k$ with $k \geq \max_{B_r} U - d$, where $A_{k, r}:= \{x \in B_r ; U(x)>k\}$, and $B_r$ is a ball of radius $r$. 
\end{lemma}

{\emph{Proof. }}
For fixed $m \geq 2$, first we show the statement for $U = \widetilde{u}_m + u_{m-2}$. 
We set $\zeta = \eta^2 \max \{u+u_{m-2}-k, 0\}$ in Proposition \ref{subsolution},  
where $u := \widetilde{u}_m$, $k$ is a real number with $k \geq \max_{B_r} (u+u_{m-2}) - d$, $\eta$ is smooth function with $\hbox{spt}\,\eta \subset B_r$, 
$0 \leq \eta \leq 1$, 
$\eta \equiv 1$ on $B_s$, $|\nabla \eta| \leq 2/(r-s)$ in $B_r \setminus B_s$, and $s = r -\sigma r \in (0, r)$, $\sigma \in (0, 1)$. Then, using the boundedness of $u, u_{m-1}, u_{m-2}$, we get 
\begin{align*}
0& \leq -\int_{A_{k, r} \cap \mathcal{S}_m(u)} \frac{u-2u_{m-1}+u_{m-2}}{h^2} \eta^2 (u+u_{m-2}-k) \,dx \\
&\quad - \frac{1}{2}\int_{A_{k, r}} (\nabla u + \nabla u_{m-2}) \cdot (2 \eta)  \nabla \eta 
(u+u_{m-2} - k) \,dx\\
&\qquad - \frac{1}{2}\int_{A_{k, r}} |\nabla u + \nabla u_{m-2}|^2 \eta^2 \,dx \\
& \leq C |A_{k, r}| + \frac{1}{2} \Bigl(\frac{1}{2} \int_{A_{k, r}} |\nabla (u+u_{m-2})|^2 \eta^2\,dx
+ 2 \int_{A_{k, r}} |\nabla \eta|^2 (u+u_{m-2} -k)^2\,dx \Bigr) \\ 
&\quad - \frac{1}{2}\int_{A_{k, r}} |\nabla u + \nabla u_{m-2}|^2 \eta^2 \,dx \\
&\leq C \Bigr[ 1 +  \frac{1}{(\sigma r)^2} \sup_{B_r} (u+u_{m-2} -k)^2 \Bigr] |A_{k, r}|
-\frac{1}{4} \int_{A_{k, s}} |\nabla (u+u_{m-2})|^2  \,dx, 
\end{align*}
where the constant $C$ depends only on $h, \mu, d, \Omega$. 

Next, we prove the same inequality for $U = -(\widetilde{u}_m + u_{m-2})$. 
Note that $-\widetilde{u}_{m}$ is a minimizer of the following functional:
\begin{align*}
&J_m^-(w) := \int_{\Omega \cap \mathcal{S}_m^-(w)}
\frac{|w +2u_{m-1} -u_{m-2}|^2}{2 h^2} \,dx 
+ \frac{1}{4} \int_{\Omega} 
|\nabla w - \nabla u_{m-2}|^2 \,dx. \\
& \qquad \text{in the set} \quad {\mathcal K}^- := \left\{ w \in H^1(\Omega); \, w = -u_0 \,\,\hbox{on}\,\,\partial \Omega \right\}, 
\end{align*}
where $\mathcal{S}_m^-(w)$ is defined to be the set $\{w<0\} \cup \{u_{m-1}>0\} \cup \{u_{m-2} >0\}$.

Now, for $w : = -\widetilde{u}_m$, we set $\varphi
:= w-\zeta \in {\mathcal K}^-$ where
$\zeta := \eta \max\{w-u_{m-2} -k, 0\}$, $k$ is a real number with $k \geq \max_{B_r} (w-u_{m-2}) - d$, and 
$\eta$ is a smooth function chosen in the same way as above. Then, by the minimality of $w$, 
\begin{align}
\label{regular_proof1}
0 &\leq J_m^-(\varphi) -J_m^-(w) \nonumber \\
&\leq  \int_{\Omega \cap \mathcal{S}_m^- (\varphi)} 
\Bigl( \frac{-2(w + 2 u_{m-1} - u_{m-2})}{2 h^2} \zeta + \frac{|\zeta|^2}{2h^2} \Bigr) \,dx \nonumber \\
&\qquad 
+\int_{\Omega} \frac{|w + 2 u_{m-1} - u_{m-2}|^2}{2 h^2} \Bigl(\chi_{\mathcal{S}_m^- (\varphi)}
- \chi_{\mathcal{S}_m^-(w)} \Bigr) \,dx  \nonumber \\
&\qquad\qquad + \frac{1}{4} \int_{\Omega} |\nabla \varphi
 - \nabla u_{m-2}|^2\,dx  - \frac{1}{4} \int_{\Omega}
|\nabla w - \nabla u_{m-2}|^2 \, dx. 
\end{align}

Note that the term in the third line is less than or equal to 
$\frac{1}{2h^2} \int_{\hbox{spt}\,\zeta} |w + 2 u_{m-1} - u_{m-2}|^2 \,dx$, since 
$\chi_{\mathcal{S}_m^- (\varphi)}- \chi_{\mathcal{S}_m^-(w)}$ is positive only for $x$ satisfying
 $0 \leq w(x) <\zeta(x)$. 
Therefore, recalling that $\hbox{spt} \,\zeta \subset A_{k, r}$, the first two terms on 
the right-hand side of (\ref{regular_proof1})
are less than or equal to $C |A_{k, r}|$, 
where $C$ is a constant depending only 
$\Omega, \mu, d, h$. Then, we continue the estimate (\ref{regular_proof1}) as follows: 
\begin{align*}
0 & \leq C |A_{k, r}| + \frac{1}{2} \int_{A_{k, r}} (1-\eta)^2 |\nabla w - \nabla u_{m-2}|^2 \,dx \\
& \quad + \frac{1}{2} \int_{A_{k, r}} (w - u_{m-2} - k)^2 |\nabla \eta|^2 \,dx 
 - \frac{1}{4} \int_{A_{k, r}} |\nabla w - \nabla u_{m-2}|^2 \, dx \\
& \leq C |A_{k, r}| + \frac{1}{2} \int_{A_{k, r}} |\nabla (w - u_{m-2})|^2 \,dx 
 + \frac{2}{(\sigma r)^2} \int_{A_{k, r}} (w - u_{m-2} - k)^2 \,dx \\
&\quad - \frac{3}{4} \int_{A_{k, s}} |\nabla (w -  u_{m-2})|^2 \, dx. 
\end{align*}
Therefore, we get 
\begin{align*}
&\int_{A_{k, s}} |\nabla (w -  u_{m-2})|^2 \, dx 
\leq C |A_{k, r}| + \theta \int_{A_{k, r}} |\nabla (w - u_{m-2})|^2 \,dx \\
&\qquad+ \frac{8}{3}\frac{1}{(\sigma r)^2} \int_{A_{k, r}} (w - u_{m-2} - k)^2 \,dx, 
\end{align*}
where $\theta = \frac{2}{3}<1$.  By Lemma V. 3.1 in \cite{[MG]}, 
we obtain
\begin{equation*}
\int_{A_{k, s}} |\nabla (w -  u_{m-2})|^2 \, dx 
\leq C |A_{k, r}| + \frac{8}{3} \frac{1}{(\sigma r)^2} \int_{A_{k, r}} (w - u_{m-2} - k)^2 \,dx, 
\end{equation*}
which is the desired estimate for $- (\widetilde{u}_m + u_{m-2})$.

{\emph{Proof of Theorem \ref{regular}}.}
Lemmas \ref{bounded} and \ref{Cacciopoli_ineq} imply 
 that $\widetilde{u}:=\widetilde{u}_{m} + u_{m-2}$ ($m \geq 2$) 
belongs to the De Giorgi class $\mathcal{B}_2(\Omega, \mu, \gamma, d)$. 
Thus, by De Giorgi's embedding theorem (\cite[Section 2.6]{[L-U]}), 
$\widetilde{u}_m + u_{m-2} \in C^{0, \widetilde{\alpha}}(\widetilde{\Omega})$ 
for some $\widetilde{\alpha} \in (0, 1)$ which is independent of $m$. We can now prove that 
$\widetilde{u}_m \in C^{0, \alpha}(\widetilde{\Omega})$ for some $\alpha \in (0, 1)$. To this end, set $\alpha := \min \{\alpha_0, \widetilde{\alpha}\}$. 
For $m=2$, by the fact that $\widetilde{u}_2 + u_0 \in C^{0, \widetilde{\alpha}}(\widetilde{\Omega}) \subset C^{0, \alpha}(\widetilde{\Omega})$, and  
the assumption $u_0 \in C^{0, \alpha_0}(\widetilde{\Omega}) \subset C^{0, \alpha}(\widetilde{\Omega})$, we see that $\widetilde{u}_2 \in
C^{0, \alpha}(\widetilde{\Omega})$. Hence, $u_2
:= \max \{\widetilde{u}, 0\}$ belongs to the same space. Now, we assume that $\widetilde{u}_k 
\in C^{0, \alpha}(\widetilde{\Omega})$ for any $k=2, ...., m-1$. Then, since $u_{m-2} \in 
C^{0, \alpha}(\widetilde{\Omega})$, and $\widetilde{u}_m + u_{m-2} \in C^{0, \widetilde{\alpha}}
(\widetilde{\Omega})
\subset C^{0, \alpha}(\widetilde{\Omega})$, 
we have $\widetilde{u}_m \in C^{0, \alpha}(\widetilde{\Omega})$.

By the above theorem, we can choose the support of test functions within the open set 
$\{ \widetilde{u}_m>0\}$, which leads to the following first variation formula for $J_m$. 

\begin{proposition}[First variation formula]
\label{first_vari_form}
Any minimizer $u$ of $J_m$ for 
$m=2,3, \dots ,M$, satisfies the following equation:
\begin{equation}
\label{first_vari}
\int_{\Omega} \left(
\frac{  u - 2u_{m-1} +u_{m-2} }{h^2}\phi
+  \nabla \frac{u+u_{m-2}}{2} \cdot \nabla \phi
\right) dx = 0
\end{equation}
for all $\phi \in C_0^\infty(\Omega\cap \{u>0\})$. 
\end{proposition}

{\emph{Proof}}
Since $\{u>0\}$ is an open set by Theorem 3.2, we can calculate 
the first variation of $J_m$ using $u+\varepsilon \phi$ with $\phi \in C_0^\infty(\Omega\cap \{u>0\})$ 
as a test function. The result then follows by noting that there exists $\varepsilon_0>0$ such that 
$\chi_{\mathcal{S}_m(u+\varepsilon \phi)} = \chi_{\mathcal{S}_m(u)}$ for  
$|\varepsilon|< \varepsilon_0$.

\section{Definition and existence of weak solutions}
\label{sec:weak solution}
In this section, we will construct weak solutions to Problem \ref{prob} in the one dimensional setting. 
First, we state the definition.

\begin{definition}[Weak solution]
\label{weak_sol}
For a given $T>0$, a \emph{weak solution} is defined as a function $u \in H^1((0, T) ; L^2(\Omega)) \cap L^\infty((0, T) ; H^1_0(\Omega))$ satisfying the following equality, for all test functions 
$\phi \in C_0^\infty(\Omega \times [0,T)\cap \{u>0\})$:
\begin{equation}
\int_0^T \int_{\Omega} \left(- u_t \phi_t
+ \nabla u \cdot \nabla \phi \right) dx dt
- \int_\Omega v_0 \phi(x,0) dx =  0. 
\end{equation}
Moreover, we require that $u \equiv 0$ is satisfied outside 
of $\{u>0\}$, and that $u(0, x) = u_0(x)$ in $\Omega$ in the sense of traces.
\end{definition}

\begin{remark}
This weak solution contains two pieces of information, namely, 
the wave equation on the positive part $\{u>0\}$, and harmonicity on the interior of the 
complement. 
If we assume the above weak solution preserves energy and has a regular free boundary,  
we can formally derive a free boundary condition solely from the definition of the weak solution. 
If we consider more general settings, such as including an adhesion term, the problem becomes 
 more complicated and requires a different notion of a weak solution. For details, see 
 Remark \ref{rem:adhesion}.   
\end{remark}

In constructing our weak solution, we carry out interpolation in time of the cut-off 
minimizers $\{ u_{m}\}$ of $J_m$, and introduce the notion of approximate weak solutions.
In particular, we define $\overline{u}^h$ and $u^h$ on $\Omega \times (0,T)$ by
\begin{eqnarray*}
\overline{u}^h(x,t) &=& u_{m}(x), \qquad m=0,1,2, \dots , M \\
u^h(x,t) &=& \displaystyle\frac{t-(m-1)h}{h}u_{m}(x)
  + \frac{mh-t}{h}u_{m-1}(x), \qquad m=1,2,3, \dots , M
\end{eqnarray*}
for $(x,t) \in \Omega \times ((m-1)h,mh]$.
These functions allow us to construct
the following approximate solution based on the first variation formula (Proposition \ref{first_vari_form}).
\begin{definition}[Approximate weak solution]
We call a sequence of functions $\{u_m\} \subset \mathcal{K}$ an
\emph{approximate weak solution} of Problem \ref{prob} if the functions $\overline{u}^h$ and $u^h$ defined above
satisfy 
\begin{eqnarray}
&&
\int_{h}^T \int_\Omega
\left(
\frac{  u^h_t(t)  - u^h_t(t-h) }{h} \phi
+
\nabla  \frac{\overline{u}^h(t)+\overline{u}^h(t-2h)}{2} \cdot \nabla \phi
\right) \, dx \, dt = 0 \quad \nonumber \\
&&
\qquad \qquad \hbox{for all}\,\,\, \phi \in C_0^\infty(\Omega \times [0,T) \cap \{u^h>0\}),
\nonumber \\
&&
\qquad \qquad
u^h \equiv 0 \quad  \mbox{\rm in}
\quad
\Omega \times (0,T) \setminus \{u^h >0\} \label{weak_approx_sol}.
\end{eqnarray}
We further require that the initial conditions $u^h(x, 0)=u_0(x)$ and 
$u^h(x, h)=u_0(x)+ h v_0(x)$ are fulfilled.
\end{definition}


If one can pass to the limit as $h \rightarrow 0$,
then the above approximate weak solutions are expected to converge to a 
weak solution defined above. In the one-dimensional setting, 
that is $\hbox{dim}\,\Omega =1$, by energy 
estimate (\ref{EE0}) in Section \ref{sec:minimizing method}, 
we obtain the following convergence result, as in \cite{[Kikuchi]}.

\begin{lemma}[Limit of approximate weak solution]
\label{limit_process}
Let $\Omega \subset \mathbb{R}$ be a bounded open interval. Then, there exists a decreasing sequence 
$\{h_j \}_{j=1}^{\infty}$ with $h_j \to 0+$ (denoted as $h$ again) and $u \in H^1((0, T) ; L^2(\Omega)) \cap L^\infty((0, T) ; 
H^1_0(\Omega))$ such that
\begin{align}
u^h_t &\rightharpoonup u_t \quad \hbox{weakly $*$ in\,\,}  L^\infty((0, T) ; L^2(\Omega)), 
\label{limit1} \\
\nabla \overline{u}^h &\rightharpoonup \nabla u \quad \hbox{weakly $*$ in\,\,}  L^\infty((0, T) ; 
L^2(\Omega)), \label{limit2} \\
u_h  & \rightrightarrows u \quad \hbox{uniformly on\,\,}  [0, T) \times \Omega \label{limit3}.
\end{align}
Moreover, $u$ is continuous on $\Omega \times (0, T)$, and satisfies the initial condition 
$u(x, 0) = u_0(x)$. 
\end{lemma}

{\emph{Proof. }} Rewriting the energy estimate (\ref{EE0}) with $\overline{u}^h$ and $ u^h$, 
we have 
\begin{equation}
\label{EE_approx}
\|u^h_t(t)\|^2_{L^2(\Omega)} + \|\nabla \overline{u}^h(t)\|_{L^2(\Omega)}^2 \leq C \quad
\hbox{for a.e.} \,\,t \in (0, T),
\end{equation}
which together with the fact that $u^h - u_0$ has zero trace on $\partial\Omega$ immediately implies (\ref{limit1}) and (\ref{limit2}). Regarding (\ref{limit3}), we first prove
the equicontinuity of the family  $\{ u^h \}$ using (\ref{EE_approx}) and the fact that, when 
$\Omega$ is an interval, for any 
$f \in H^1_0(\Omega)$, we have
$$
\|f\|_{L^\infty(\Omega)} \leq C \|f\|_{L^2(\Omega)}^{1/2} \|f' \|_{L^2(\Omega)}^{1/2}.
$$
Indeed, for any $t,s \in [0,T)$
\begin{eqnarray*}
\| u^h(t) - u^h(s) \|_{L^{\infty}(\Omega)}^2 &\leq& C \| u^h(t)-u^h(s) \|_{L^2(\Omega)} \| \nabla u^h(t) - \nabla u^h(s) \|_{L^2(\Omega)} \\
&\leq& C \int_s^t \| u_t^h(\tau) \|_{L^2(\Omega)} \, d\tau \\
&\leq& C |t-s|,
\end{eqnarray*}
and thus
\begin{eqnarray*}
| u^h(b,t)-u^h(a,s)| &\leq& | u^h(b,t)-u^h(a,t)| + | u^h(a,t)-u^h(a,s)| \\
&=& \left| \int_a^b \frac{\partial}{\partial x} u^h(\xi, t) \, d\xi \right| + | u^h(a,t)-u^h(a,s)| \\
&\leq& \| \nabla u^h \|_{L^{\infty}((0,T);L^2(\Omega))} |b-a|^{1/2} + C| t-s |^{1/2} .
\end{eqnarray*}

Moreover, the uniform boundedness of the family $\{u^h\}$ follows as a by-product. 
Therefore, invoking the Ascoli-Arzel\`a theorem concludes the proof of (\ref{limit3}).

The following lemma is needed to prove the existence of weak solutions. 
\begin{lemma}
\label{limit_remark}
Under the assumption of Lemma \ref{limit_process}, 
define $\overline{w}^h(x, t) :=0$ if $t \in (0, h]$, and $\overline{w}^h(x, t):=\overline{u}^h(t-2h)$ when $t \in(h, T)$. Then, 
$$
\nabla \overline{w}^h \rightharpoonup \nabla u \quad \hbox{weakly $*$ in\,\,}  
L^\infty((0, T) ; L^2(\Omega)). 
$$
\end{lemma}

{\emph{Proof. }}
In the following argument, we omit the space variable $x$ for simplicity.  
We fix $U \in L^1((0, T) ; L^2(\Omega))$ and extend it by zero outside of $(0, T)$. The extended 
function, denoted again by $U$, belongs to $L^1((-\infty, \infty) ; L^2(\Omega))$. 
We calculate as follows: 
\begin{align}
\label{limit_rem_proof1}
&\Bigl |\int_{0}^T \langle \nabla \overline{w}^h(t), U(t) \rangle_{L^2(\Omega)} \,dt 
-\int_{0}^T \langle \nabla u(t), U(t) \rangle_{L^2(\Omega)} \,dt \Bigr|  \nonumber \\
&= \Bigl |\int_{-h}^{T-2h} \langle \nabla \overline{u}^h(t), U(t+2h) \rangle \,dt 
-\int_{0}^T \langle \nabla u(t), U(t) \rangle \,dt \Bigr| \nonumber \\
&\leq \Bigl |\int_{0}^{T-2h} \langle \nabla \overline{u}^h(t), U(t+2h) - U(t) \rangle \,dt \Bigr| \nonumber  \\
& \quad +\Bigl |\int_{0}^{T}  \langle \nabla \overline{u}^h(t)-\nabla u(t), U(t) \rangle \,dt\Bigr| 
+\Bigl |\int_{T-2h}^{T}  \langle \nabla \overline{u}^h(t)-\nabla u(t), U(t) \rangle \,dt\Bigr| \nonumber  \\
& \quad +\Bigl| \int_{-h}^{0} \langle \nabla \overline{u}^h(t), U(t+2h) \rangle \,dt\Bigr| 
+ \Bigl| \int_{T-2h}^{T}  \langle \nabla u(t), U(t) \rangle \,dt\Bigr| \nonumber \\
&\leq C \int_{-\infty}^{\infty} \|U(t+2h) - U(t)\|_{L^2(\Omega)} \,dt
+\Bigl |\int_{0}^{T}  \langle \nabla \overline{u}^h(t)-\nabla u(t), U(t) \rangle \,dt\Bigr| \nonumber \\
&\qquad + C \int_{T-2h}^{T} \|U(t) \|_{L^2(\Omega)} \,dt 
+ C \int_{h}^{2h} \|U(t) \|_{L^2(\Omega)} \,dt, 
\end{align}
where the constant $C$ is independent of $h$. Letting $h \to 0+$, the second term 
converges to $0$ by (\ref{limit2}), and the remaining terms vanish thanks to 
the integrability of $U$.

We now arrive at the following theorem:

\begin{theorem}[Existence weak solutions to Problem \ref{prob}]
Let $\Omega$ be a bounded domain in $\mathbb{R}$. Then Problem \ref{prob}
 has a weak solution 
in the sense of Definition \ref{weak_sol}. 
\end{theorem}

{\emph{Proof. }}
The proof is similar to that in \cite{[KaO]} and \cite{[G-S]}. 
Without loss of generality, we can consider 
$\Omega = (0, 1)$. By the definition of an approximate
weak solution (\ref{weak_approx_sol}), we have 

\begin{eqnarray}
&&
\int_{h}^T \int_\Omega
\left(
\frac{u^h_t(t)  - u^h_t(t-h)}{h} \varphi
+
\nabla  \frac{\overline{u}^h(t)+\overline{u}^h(t-2h)}{2} \cdot \nabla \varphi
\right) dx dt = 0 \quad \nonumber \\
&&
\qquad \qquad \forall \varphi \in \mathcal{C}(\overline{u}^h) := C_0^\infty(\Omega \times [0,T) \cap 
\{\overline{u}^h>0\}),
\nonumber \\
&& \label{eqh1}
\qquad \qquad
u^h \equiv 0 \quad  \mbox{\rm in}
\quad
\Omega \times (0,T) \setminus \{u^h >0\} \label{weak_sol_1}.
\end{eqnarray}
We fix $\psi \in \mathcal{C}(u)$, where $u$ is obtained in Lemma \ref{limit_process}. 
Since $u$ is continuous on $\Omega \times (0, T)$, 
there exists $\eta > 0$ such that $u \geq \eta$ on $\hbox{spt}\,\psi$. 
By Lemma \ref{limit_process}, the subsequence $\{u^h\}$ converges to $u$ uniformly, and there exists
$h_0 >0$ such that 
$$
\max_{(x, t) \in \Omega \times (0, T)} |u^h(x, t) -u(x, t)| \leq \frac{\eta}{2} \quad \hbox{for all}\,\,
\, h < h_0.
$$
Therefore, we have $u^h \geq u -|u^h -u| \geq \eta/2$ on $\hbox{spt}\,\, \psi$ for any $h \in
 (0, h_0)$. Note that $\overline{u}^h(x, t) = u^h(x, kh)$ for any $t \in ((k-1)h, kh]$, and 
$\overline{u}^h
 \geq \eta/2 >0$ on spt $\psi$ for any $h \in (0, h_0)$. This implies that (\ref{weak_sol_1}) holds 
 for any test function $\varphi \in \mathcal{C}(u)$ whenever $h<h_0$.  
 The time-discretized term can be rearranged as
\begin{multline*}
 \int_{h}^T \frac{u^h_t(t)  - u^h_t(t-h)}{h} \varphi (t) \, dt 
=  \int_h^T u_t^h(t) \frac{\varphi (t) - \varphi(t+h)}{h} \, dt \\
 - \frac{1}{h} \int_0^h u_t^h(t) \varphi(t+h) \, dt + \frac{1}{h}\int_{T-h}^T u_t^h(t) \varphi (t+h) \, dt .
 \end{multline*}
Hence, using Lemma \ref{limit_process} and Lemma \ref{limit_remark}, and passing to $h \to 0+$ in \eqref{eqh1}, we obtain
$$
 \int_0^T \int_{\Omega} \left(- u_t \varphi_t
+ \nabla u \cdot \nabla \varphi \right) dx dt
- \int_\Omega v_0 \varphi(x,0) dx =  0 \qquad \forall \varphi \in \mathcal{C}(u),
$$
which was our goal. 

\begin{remark}
Note that the fact that (\ref{weak_sol_1}) holds for any test function $\varphi$ compactly supported
in the support of the limit function $u$, is a consequence of the uniform convergence of the approximating 
sequence $\{ u^h\}$, which was in turn obtained thanks to Ascoli-Arzel\`a theorem through an imbedding argument.
However, this imbedding is available only in spatial dimension 1, restricting our existence result. 
We are not aware of any method to extend
this limit passage relying on imbedding theorems, except for the recent results in \cite{[BNO], [BLNO]},
where a different definition of weak solution through the subsolution property allows proving existence in arbitrary dimension.
Nevertheless, the good point of our strategy, common to \cite{[BNO], [BLNO]}, is that a uniform 
estimate is available for the approximations and hence a unique limit function can be identified. 
The open question is how much can be said about this limit.
On the other hand, the energy-preserving property of our approximation scheme could be exploited not only numerically 
in building robust numerical algorithms, but also theoretically in proving uniqueness
of solutions for certain wave-types problems. The existing analyses using variational approach were not successful 
in proving uniqueness, and hence investigating the new scheme from this viewpoint presents an interesting direction
for future research. 
\end{remark}

\begin{remark}
\label{rem:adhesion}
As a remark concluding this section, we present formal ideas on a possible approach to deal with the free boundary condition. 
For this purpose we generalize Problem \ref{prob} to the setting of a hyperbolic free boundary problem with an adhesion term, which delineates the role of the free boundary condition and is of interest in applications.
We calculate the first variation of the following action integral: 
$$
J(u) := \int_{0}^T \int_{\Omega} \Bigl( (u_t)^2 \chi_{\{u>0\}} -|\nabla u|^2 -Q^2 \chi_{\{u>0\}}\Bigr)\,dx\,dt
$$
where $Q$ is a constant which expresses the adhesion force. 
When energy is conserved, i.e., when the function $u$ does not change its value from positive to zero as time passes, we can, under appropriate assumptions, 
calculate the first variation, as well as the inner variation, of the functional $J$. 
However, if energy is not conserved, we can calculate neither the first variation nor 
the inner variation due to the presence of the $Q^2$-term containing the characteristic function. 
To overcome this difficulty, we consider a smoothing of the characteristic function within the adhesion term by a 
function $B_{\varepsilon}$ defined by $B_{\varepsilon}(u)= \int_{-1}^u \beta_{\varepsilon}(s) \, ds$, where 
$\beta_{\varepsilon}(u):=\frac{1}{\varepsilon} \beta(\frac{u}{\varepsilon})$, and
 $\beta : \mathbb{R} \to [0, 1]$ is a smooth function satisfying $\beta = 0$ outside $[-1, 1]$, $\int_{\mathbb{R}} \beta(s)\,ds =1$, and $B(0)=\frac{1}{2}$. After smoothing, we can calculate the first variation to obtain an expression for the following problem:
\begin{equation*}
(\hbox{P}_\varepsilon)
\begin{cases}
\chi_{\overline{\{u>0\}}}\, u_{tt} &= \Delta u - \frac{1}{2} Q^2 \beta_{\varepsilon} (u) 
\quad \hbox{in} \,\,\Omega \times (0, T), \\
u(x, 0) &= u_0(x) \quad \hbox{in} \,\, \Omega, \\
u_t(x, 0) &= v_0(x) \quad \hbox{in} \,\, \Omega,\\
u(x, t) |_{\partial \Omega} &= \psi(x, t) \,\,\hbox{with}\,\,\psi(x, 0) = u_0 \,\,\,\hbox{on}\,\,
 \partial \Omega,
\end{cases}
\end{equation*}
where $u_0, v_0$ are the same as in Problem \ref{prob}, and $\psi$ is a given function. Now, we set the following hypotheses: 
\begin{description}
\item[(H1)] The existence of a solution $u_\varepsilon$ to $(\hbox{P}_\varepsilon)$. 
\item[(H2)] The existence of a function $u : \Omega \times (0, T) \to \mathbb{R}$ such that 
$u_\varepsilon \to u$ in an appropriate topology as $\varepsilon \downarrow 0$ and such that the
 following holds:
	\begin{description}
	\item[(H2.1)] $u_{tt} -\Delta u=0$ in $\Omega \times (0, T) \cap \{u>0\}$. 
	\item[(H2.2)] The free boundary $\partial \{u>0\}$ is regular, $\mathscr{H}^{N}(\mathcal{D} \cap 
	\partial \{u>0\})<\infty$ for any $\mathcal{D} \subset \subset \Omega \times (0, T) \subset 
	\mathbb{R}^N \times (0, T)$, and 
	$|D u| \not= 0$ on $\Omega \times (0, T) \cap \partial \{u>0\}$. Here,  
	$Du = (u_{x_1}, ..., u_{x_N}, u_t)$, and $\mathscr{H}^N$ is the $N$-dimensional Hausdorff measure. 
	\item[(H2.3)] $u$ is a subsolution in the following sense:
	$$
	\int_{0}^T \int_{\Omega} \Bigl( \chi_{\overline{\{u>0\}}}\,\, u_{tt}\, \zeta
	+\nabla u \cdot \nabla \zeta \Bigl) \, dx\,dt \leq 0
	$$
	for arbitrary nonnegative $\zeta \in C_0^{\infty}(\Omega \times (0, T))$. 
	\end{description}
\end{description}

Starting from $(\hbox{P}_\varepsilon)$, and employing (H1), (H2.1) and (H2.2), 
we can show that the limit function $u$ satisfies the following free boundary condition \cite{[O4]}:
\begin{equation}
\label{free_boundary_cond}
|\nabla u|^2 - u_t^2 =Q^2 \quad \hbox{on} \,\,\,\Omega \times (0, T) \cap \partial \{u>0\}. 
\end{equation}
Now, for any $\mathcal{D} \subset \subset \Omega \times (0, T)$, we define a linear 
functional $f$ on $C_0^{\infty}(\mathcal{D})$ corresponding to 
$\Delta u - \chi_{\overline{\{u>0\}}} \,u_{tt}$ 
as follows:
$$
f (\zeta) := - \int_{\mathcal{D}} \Bigl( \chi_{\overline{\{u>0\}}}\,\, u_{tt}\, \zeta
+\nabla u \cdot \nabla \zeta \Bigr)\, dx\,dt.
$$
Since $f$ is a positive linear functional on $C_0^{\infty}(\mathcal{D})$ by (H2.3), 
$f$ can be extended to a positive linear functional on $C_0(\mathcal{D})$. 
Riesz's representation theorem asserts there exists a unique 
positive Radon measure $\mu_f$ on $\mathcal{D}$ such that 
$ f(\zeta) = \int_{\mathcal{D}} \zeta \,d\mu_f $. 
In this sense, we can say that $\Delta u - \chi_{\overline{\{u>0\}}}\,\, u_{tt}$ is 
a positive Radon measure on $\mathcal{D}$. 

On the other hand, we can calculate the value of $f(\zeta)$ from
 (\ref{free_boundary_cond}). By splitting the integral domain into four parts, 
$\mathcal{D} \cap \{u>0\}$, $\mathcal{D} \cap \partial \{u>0\}$, 
$\mathcal{D} \cap \partial \{u=0\}^{\circ}$, $\mathcal{D} \cap \{u=0\}^{\circ}$, noting that 
all terms vanish except for $\mathcal{D} \cap \partial \{u>0\}$ by the integration by parts, and 
$\chi_{\overline{\{u>0\}}}=1$ on $\partial \{u>0\}$, 
we get
\begin{align*}
&\int_{\mathcal{D}} \left( \chi_{\overline{\{u>0\}}}\,\, u_{tt} \zeta
+\nabla u \cdot \nabla \zeta \right) \, dx\,dt 
= \int_{\mathcal{D}} \left( -(\chi_{\overline{\{u>0\}}}\,\,\zeta)_t u_{t} 
+\nabla u \cdot \nabla \zeta \right) \, dx\,dt \\
&= \int_{\mathcal{D} \cap \partial \{u>0\}} -\Bigl(\chi_{\overline{\{u>0\}}}\,\,\zeta \cdot u_{t}
\cdot \frac{-u_t}{|Du|} 
+\zeta \nabla u \cdot \frac{-\nabla u}{|Du|}\Bigr) \, d\mathscr{H}^N\\
&= \int_{\mathcal{D} \cap \partial \{u>0\}}  
\frac{u_t^2-|\nabla u|^2}{|Du|}\,\, \zeta  \, d\mathscr{H}^N 
=\int_{\mathcal{D} \cap \partial \{u>0\}}  
\frac{-Q^2}{|Du|}\,\, \zeta  \, d\mathscr{H}^N  \quad (\hbox{by} \,\,\,(\ref{free_boundary_cond})). 
\end{align*}
From the above and the definition of $f$,  we observe that 
\begin{equation}
\label{measure_represen3}
\mu_f = \frac{Q^2}{|Du|} \mathscr{H}^N \lfloor \partial \{u>0\}.
\end{equation}
In this sense, the positive Radon measure $\Delta u - \chi_{\overline{\{u>0\}}}\,\, u_{tt}$ has its 
support in the free boundary $\partial \{u>0\}$. Formally, we can rewrite (\ref{measure_represen3}) as follows:
\begin{equation}
\label{hfbp_eq}
\chi_{\overline{\{u>0\}}} \,u_{tt} - \Delta u  = -\frac{Q^2}{|Du|} \mathscr{H}^N \lfloor \partial \{u>0\}. 
\end{equation}

Summarizing the above, starting from the smoothed problem ($P_\varepsilon$), 
under the hypotheses (H1)-(H2), we formally derive a hyperbolic degenerate equation with 
adhesion force.  This equation (\ref{hfbp_eq}) includes all information about the 
hyperbolic free boundary problem, that is, the wave equation $u_{tt} -\Delta u = 0$ in the set $\{u>0\}$, 
the free boundary condition $|\nabla u|^2 - u_t^2 = Q^2$ on
$\Omega \times (0, T) \cap \partial \{u>0\}$, and 
the Laplace equation $\Delta u= 0$ in the set $\{u<0\}$ a.e. $t \in (0, T)$. 
\end{remark}

\section{Numerical Results}
\label{sec:numerics}
\subsection{Numerical analysis of the one-dimensional problem}
In this section, we present several numerical results for the equation \eqref{prob} obtained by the minimization of the Crank-Nicolson type functional 
\begin{align*}
J_m(u) := &\int_{\Omega  \cap (\{u>0\} \cup \{u_{m-1}>0\} \cup \{u_{m-2}>0\})}
\frac{|u -2u_{m-1} +u_{m-2}|^2}{2 h^2} \,dx  \\
&\qquad+ \frac{1}{4} \int_{\Omega} 
|\nabla u + \nabla u_{m-2}|^2 \,dx. 
\end{align*}
and compare them with results due to the original discrete Morse flow method of \cite{[O3]}, which uses the functional 
$$
\widetilde{J}_m(u) := \int_{\Omega \cap (\{u>0\} \cup \{u_{m-1}>0\} )}
\frac{|u -2u_{m-1} +u_{m-2}|^2}{2 h^2} \,dx + \frac{1}{2} \int_{\Omega} 
|\nabla u|^2 \,dx. 
$$
In the numerical calculation, we simply use the functional $I_m$ without the restriction of the integration domain and the corresponding functional $\widetilde{I}_m$ for the original discrete Morse flow method. Subsequently, 
for a minimizer $\widetilde{u}_m, \, m \geq 2$, of $I_m$ or $\widetilde{I}_m$, we define 
$$
u_m:=\max\{\widetilde{u}_m,0\}. 
$$
We regard $u_m$ as a numerical solution at time level $t=mh$. 
The minimization problems are discretized by the finite element method, where 
the approximate minimizer is a continuous function over the domain and piece-wise linear over each element. 

In the one-dimensional case, equation \eqref{prob} has been employed in describing the dynamics of a string hitting a plane with zero reflection constant. In two dimensions, the graph of the solution may be considered as representing a soap film touching a water surface. Another important application of this model is the volume constrained problem describing the motion of scalar droplets over a flat surface (see, e.g., \cite{[G-S]}, 
and Section \ref{sec:higherdim}). 

\begin{figure}[htp]
\begin{center}
	\includegraphics[width=0.30\textwidth]{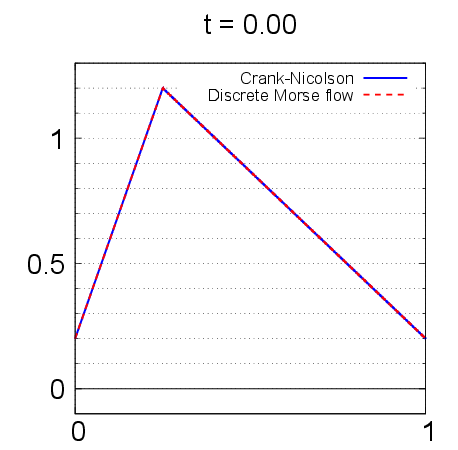}
	\includegraphics[width=0.49\textwidth]{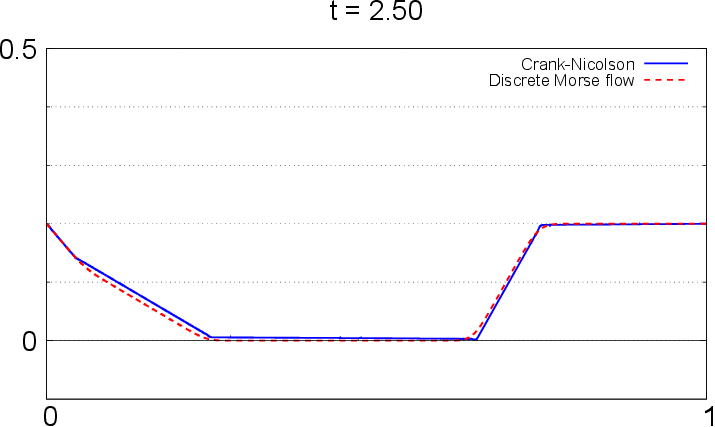}
	\\
	\includegraphics[width=0.49\textwidth]{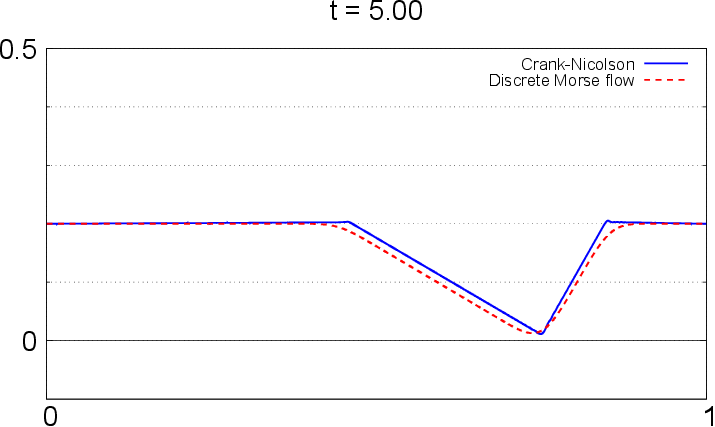}
	\includegraphics[width=0.49\textwidth]{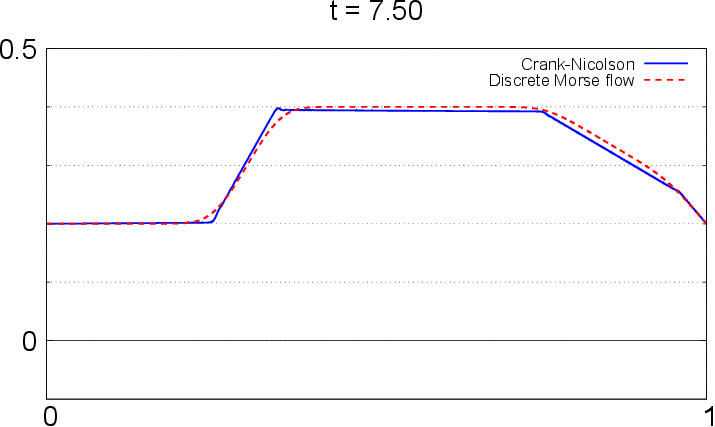}
	\caption{Numerical solution at four distinct times for the Crank-Nicolson method (blue) and the original discrete Morse flow method (red). The time step size is $h = 1.0 \times 10^{-4}$ and the
	spacial mesh size is $\Delta x = h$. }
	\label{fig:u}
  \end{center}
\end{figure}

Having in mind the model of a string hitting an obstacle, let us first consider problem \eqref{prob} in the open interval $\Omega=(0,1)$, with the initial condition 
\begin{equation*}
u_0(x) := \left\{
\begin{aligned}
&4x+0.2	&&\text{ if } 0\leq x <1/4,\\
&-\frac{4}{3}(x-1)+0.2	&&\text{ otherwise },
\end{aligned}
\right.
\end{equation*}
and $v_0\equiv 0$. 
Figure \ref{fig:u} shows the behavior of the numerical solution for both methods.
For the Crank-Nicolson method, the corners in the graph of the solution are kept, 
even as time progresses.
This is not the case for the discrete Morse flow method, where corners are smoothed.

\begin{figure}[htp]
\begin{center}
	\includegraphics[width=0.85\textwidth]{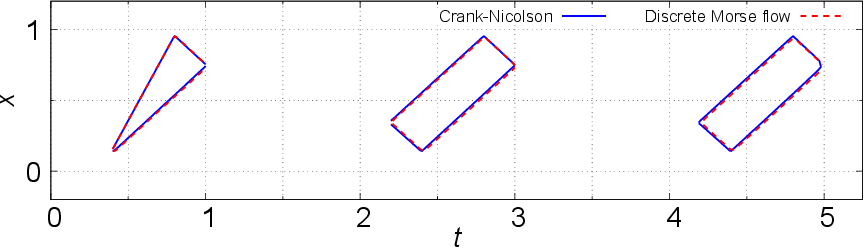}
	\vspace{-0.cm}
	\caption{Free boundary corresponding to the motion in Figure \ref{fig:u}. 
The curves are obtained by plotting the boundary of the set $\{(t,x);u(x,t) < \varepsilon\}$ for a small $\varepsilon >0$.}
	\label{fig:fb}
\end{center}
\end{figure}

\vspace{0.3cm}
Figure \ref{fig:fb} shows that the free boundary condition \eqref{eq:fbc0} is satisfied when the string peels off the obstacle.

\begin{figure}[htp]
\begin{center}
	\includegraphics[width=0.7\textwidth]{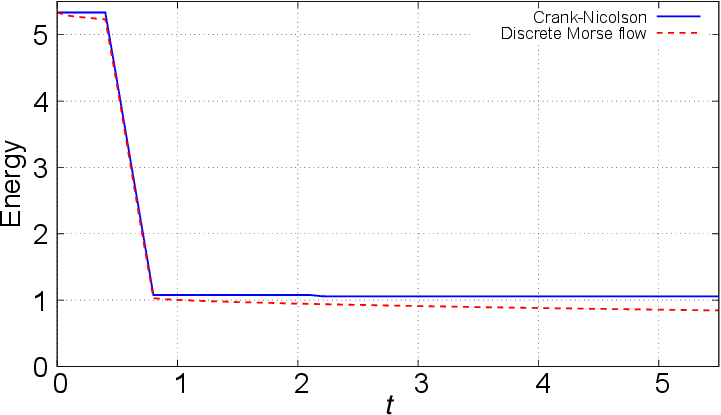}
	\caption{Evolution of the energy of the numerical solution for both methods.}
	\label{fig:ene}
\end{center}
\end{figure}

\vspace{0.3cm}
Figure \ref{fig:ene} shows that the energy is lost when the string touches the obstacle, while the energy is preserved before and after the contact of the string with the obstacle. For the sake of comparison, we note that the energy of the solution obtained by discrete Morse flow decays even during the non-contact stage. \\

\begin{figure}[htp]
\begin{center}
	\includegraphics[width=0.8\textwidth]{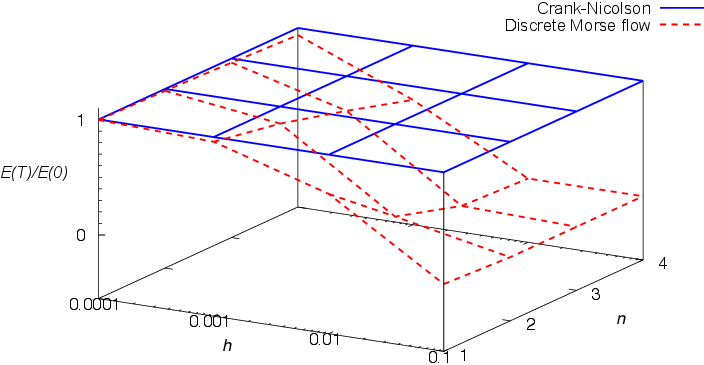}
	\vspace{10pt}
	\caption{Comparison of energy decay tendency for both methods using the initial data 
	$u_0 = \sin(2n\pi x)$ and $v_0 \equiv 0$. Here, $\Delta x = h$ is used. }
	\label{fig:ene_lost}
\end{center}
\end{figure}

\vspace{0.3cm}
To test the energy decay tendency of both methods, 
we solved the problem without free boundary with the initial condition $u_0 = \sin(2n\pi x)$, and $v_0 \equiv 0$. 
It was found that, for the original discrete Morse flow, energy decay becomes prominent with decreasing time resolution and increasing wave frequency. On the other hand, as can be observed in Figure \ref{fig:ene_lost}, the Crank-Nicolson method preserves energy independent of the time resolution and wave frequency.

\vspace{0.3cm}
Although the Crank-Nicolson method displays excellent energy-preserving properties, it appears to include an incorrect phase-shift, as is the case with the original discrete Morse flow.
We summarize the features of both methods in Table \ref{tab:feature}.

\begin{table}[h]
	\centering
	\begin{tabular}{c|cc}
		\toprule
		& C-N & DMF  \\ \hline \hline 
		energy         & conserved & decays \\ \hline 
		free boundary condition & holds & holds \\ \hline 
		high harmonic wave				   & preserved & decays \\ \hline	
		including constraints & possible & possible  \\ \hline
		phase shift   & occurs & occurs  \\
		\bottomrule \\
	\end{tabular}
	\caption{Main features of the two methods compared in this section.}
	\label{tab:feature}
\end{table}

\subsection{Higher dimensions and more general problems}
\label{sec:higherdim}
In this section, we investigate the energy preservation properties of the proposed scheme in the two dimensional setting. In particular, the functional (\ref{eq:ECfunctional}) is used to approximate a solution of the wave equation
with initial conditions $u_0(x, y) = \sin (\pi x) \sin (\pi y), v_0(x, y) = 0$ and Dirichlet zero boundary condition, where the domain $\Omega=(0,1)\times(0,1)$. 

The functional value corresponding to a given function $u$ is approximated using $P_1$ finite elements,
and the functional minimization is performed using a steepest descent algorithm. Here $\Omega$ has been uniformly partitioned into $N=5684$ elements. 

\begin{figure}[htp]
\begin{center}	
	\includegraphics[width=0.85\textwidth]{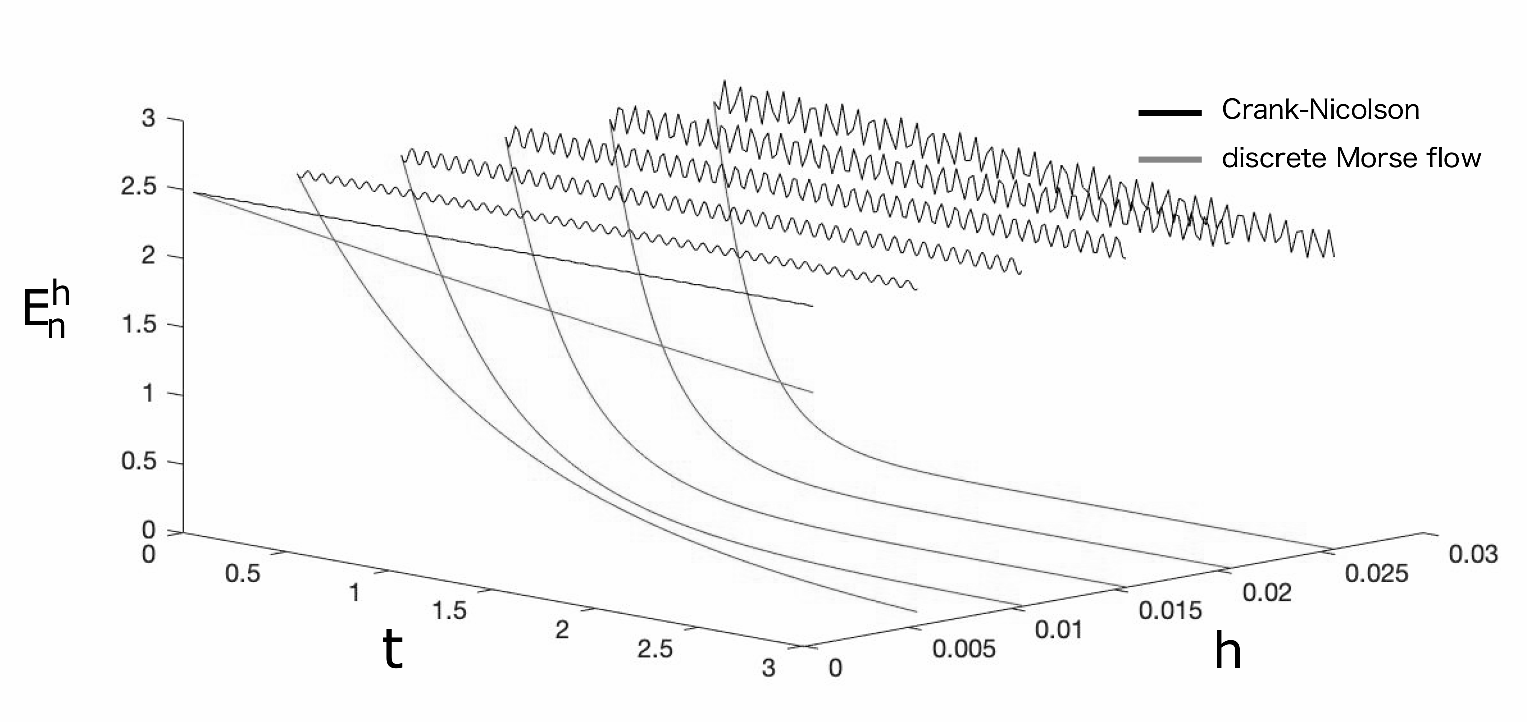}
	\caption{Comparison of the Crank-Nicolson scheme with the original discrete Morse flow for a 2-dimensional problem.}  
	\label{fig:comp}
\end{center}
\end{figure}

Using several different values of the time step $h,$ we compared the energy of the numerical solution obtained using the Crank-Nicolson scheme with that obtained from the standard discrete Morse flow. 
The total energy is computed using the finite element method on the functional:
\begin{align}
\mathcal{E}^h_n(u)=\int_{\Omega} \Bigl( \frac{1}{2}\left|\frac{u-u_{n-1}}{h}\right|^2+\frac{|\nabla u|^2}{2} \Bigr) \,dx.
\end{align}
The results are shown in Figure \ref{fig:comp}, where the time steps were $h=0.0005+0.005\times k$, $k = 0,1,2,3,4,5$. Our results confirm the energy preservation properties of the proposed scheme. 

\vspace{0.3cm}
We have also used the proposed method to investigate the numerical solution of a more complicated model equation describing droplet motions. 

\begin{figure}[htp]
 \begin{center}
	{\includegraphics[bb=0 0 500 600,scale=0.67]{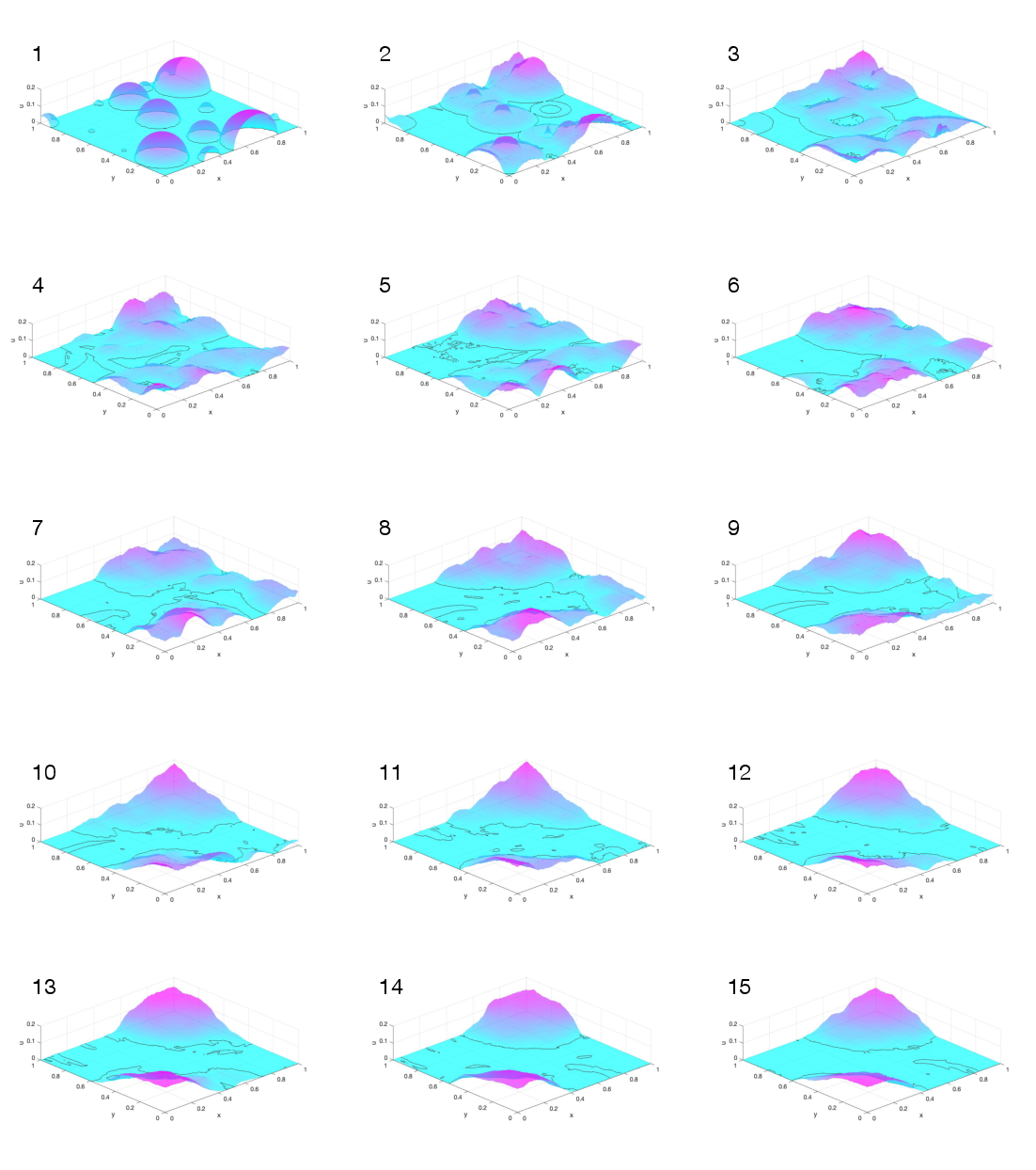}}
	\caption{Crank-Nicolson type minimizing movement approximation of droplet motion, with the free boundary illustrated as the black curves.. 
	Time is designated by the integer values within the figure, so that the initial condition corresponds 
	to number 1 and all graphs are plotted at equal time intervals, except for the last one showing the stationary state reached after sufficiently long time.} 
	\label{fig:drops}
\end{center}
\end{figure}
 
\vspace{0.3cm}
The target equations correspond to volume constrained formulations of the original problem. In particular, volume and non-negativity constraints are added to the functionals 
by means of indicator functions:
\begin{multline*}
 \mathcal{J}^i_{m}(u) = \int_{\Omega} \left(\frac{|u-2u^i_{m-1}+u^i_{m-2}|^2}{2 h^2} + 
 \frac{1}{4}  |\nabla u^i + \nabla u^i_{m-2}|^2 + \right)dx \\ + \Psi_1(u) + \Psi^{i, m}_2(u),
\end{multline*}
 where each indicator function is defined as follows:
\begin{align}
 \Psi_{1}(u) &= {\begin{cases}
 	{0,} & \text{if } {u(x) \geq 0 \text{ for } \mathcal{L}^{N}\text{-a.e. } x \in \Omega}  \\
 					{ \infty,} & {\text{otherwise}}  \\
			   	\end{cases}}\notag , \\
	\Psi_{2}^{i,m}(u) &= {\begin{cases}
  					{0,} & \text{if } {\int_{\Omega} u(x)dx = V^i_m}  \\
	    			{ \infty,} & {\text{otherwise.}}  \\
			   		\end{cases}}\notag
\end{align}
 Here $V^i_m$ denotes the volume of droplet $i$ at time step $m$. 
 
 By minimizing functionals $\mathcal{J}_m^i$ for each droplet, we are able to compute approximate solutions to the volume constrained problem. The results are shown in Figure \ref{fig:drops}. For each $i$, the initial condition is  prescribed as a collection of spherical caps, and we observe the droplets oscillate while coalescing into larger groups.

\section{Conclusions}
\label{sec:conclusions}
We have shown the existence of  weak solutions to a hyperbolic free boundary problem by minimizing a Crank-Nicolson type functional in the one-dimensional setting. This new functional was shown to preserve the energy correctly both on continuous and discrete levels, which is of significance in both theoretical analysis and numerical simulations.
Future tasks include extending our result  to higher dimensions and to developing computational methods for investigating the numerical properties of the free boundary problem.

\section*{Acknowledgments} 
E. Ginder would like to acknowledge the support of JSPS Kakenhi Grant number 17K14229.
The research of the fourth author was partially funded by a joint research project with YKK Corporation.
The research of the last author was supported by JSPS Kakenhi Grant numbers 19K03634 and 18H05481.
We would also like to thank the anonymous referees for their constructive comments.

\end{document}